\documentclass[11pt]{article} 
\usepackage{amssymb,amsfonts,amsmath,latexsym,epsf,tikz,url}

\newtheorem{theorem}{Theorem}[section]

\newtheorem{corollary}[theorem]{Corollary}

\newcommand{\qed}{\hfill $\square$\medskip}
\begin{document}

%Investigation of domination polynomials and of some families of graphs 
\title{Counting  independent dominating sets in linear polymers}

\author{Somayeh Jahari
\and
Saeid Alikhani$^{}$\footnote{Corresponding author}}

\date{}

\maketitle

\begin{center}
 Department of Mathematics, Yazd University, 89195-741, Yazd, Iran\\
%Supported by the 

{\tt s.jahari@gmail.com ~~~~ alikhani@yazd.ac.ir}
\end{center}

%%%%%%%%%%%%%%ABSTRACT%%%%%%%%%%%%%%%%%%%%%%%%%%%%%%%%%%%%%%%%%%%%%%%%%%%%%%%%%%%%

\begin{abstract}
 
A vertex subset $W\subseteq  V$ of the graph $G = (V,E)$ is an independent dominating set if every vertex in $V\setminus W$ is adjacent to at least one vertex in $W$ and the vertices of $W$ are pairwise non-adjacent. 
We enumerate independent dominating sets  in several classes
of graphs made by a linear or cyclic concatenation of basic building blocks.  Explicit recurrences are derived for  the number of  independent dominating sets
 of 
these kind of graphs. Generating functions for the
number of independent  dominating sets of  triangular and squares cacti chain are also computed.
\end{abstract}

\noindent{\bf Keywords:} Independent dominating sets, generating function, Cactus graphs.

\medskip
\noindent{\bf AMS Subj. Class.:} 05C30; 05C69

%%%%%%%%%%%%%%%%%%%%%%%%%%%%%%%%%%%%%%%%%%%%%%%%%%%%%%%%%%%%%%%%%%%%%%%%%%%%%%%%%
%%%%%%%%%%%%%%%%%%%%%%%%%%%%%%%%%%%%%%%%%%%%%%%%%%%%%%%%%%%%%%%%%%%%%%%%%%%%%%%%%
\section{Introduction}
%%%%%%%%%%%%%%%%%%%%%%%%%%%%%%%%%%%%%%%%%%%%%%%%%%%%%%%%%%%%%%%%%%%%%%%%%%%%%%%%%
 Let $G=(V(G),E(G))$ be a  simple of finite orders graph, i.e., graphs are undirected with no loops or
parallel edges and with finite number of vertices.   
For any vertex $v\in V(G)$, the {open neighborhood} of $v$ is the
set $N(v)=\{u \in V (G) | uv\in E(G)\}$ and the {closed neighborhood} of $v$
is the set $N[v]=N(v)\cup \{v\}$. For a set $S\subseteq V(G)$, the open
neighborhood of $S$ is $N(S)=\bigcup_{v\in S} N(v)$ and the closed neighborhood of $S$
is $N[S]=N(S)\cup S$. A non-empty set $S\subseteq V(G)$ is a {dominating set}, if $N[S]=V$, or equivalently,
every vertex not in $S$ is adjacent to at least one vertex in $S$, and $S$ is a {total dominating set},  if every vertex of $V$ is adjacent to some vertices of $S$. The {domination number} ({total domination number}) of  the graph $G$, denoted by $\gamma(G)$ ($\gamma_t(G)$), is  the minimum cardinality of all dominating sets (total dominating sets) of  $G$.  
An {independent set} in a graph $G$ is a set of pairwise non-adjacent vertices. A maximum independent set in $G$ is a largest independent set and its size is called {independence number} of $G$ and is denoted by $\alpha(G)$.   An independent dominating set of $G$ is a set that is both dominating and independent in $G$. Equivalently, an independent dominating set is a maximal independent set. The independence domination number of $G$, denoted by $\gamma_i(G)$, is the minimum size of all independent dominating sets of $G$.  %In particular,  \[\gamma(G)\leq \gamma_i(G)\leq \alpha(G).\]
For a detailed treatment of domination theory, the reader is referred to~\cite{domination}.

Most of the papers published so far deal with structural aspects of domination, trying to determine exact expressions for $\gamma(G)$ or some upper and/or lower bounds for it. The enumerative side of the problem has considered in last decade and the domination polynomial of a graph $G$, denoted by $D(G,x)$, which is the generating function for the number of dominating sets, has studied well in literature. A tribonacci-like recurrence for the number of dominating sets of paths, cycles and non $P_4$-free graphs with a recurrence for the number of dominating sets for general graphs, which is dependent to another polynomial, provided \cite{TMU}. Also it is interesting that   the number of dominating sets, i.e., $D(G,1)$ is odd and more precisely, 
$D(G,r)$ is odd for all odd integer $r$ (see \cite{GCOM,Brouwer}).

 After counting dominating sets,  some authors studied the number of another kind of dominating sets \cite{outer,doslic}, especially the number of total dominating sets  and independent  dominating sets has studied well, see e.g. \cite{dod1,dod2}.

 Note that  counting the number of dominating sets is \# P-complete, even in restricted graph classes such as, split graphs and bipartite chordal graphs \cite{15}. So it is natural to look to the classes with specific constructions to obtain the number of their 
 dominating sets. In this paper we consider graphs with simple connectivity patterns, for example cacti.   
 
 The goal of this  paper is to further contribute to the corpus of knowledge about the enumerative aspects of independent domination by investigating the number of
 independent  dominating sets in some  classes of simple linear polymers. 
  A linear polymer is a long continuous chain of carbon-carbon bonds with the remaining two valence bonds attached primarily to hydrogen or another relatively small hydrocarbon moiety. We consider simple linear polymers which  
  called cactus chains. Cactus graphs were first known as Husimi trees; they appeared in the scientific  literature in papers by Husimi and Riddell concerned with cluster integrals in the theory of condensation in statistical mechanics \cite{10,14,18}. In the meantime, they also found applications in chemistry \cite{13,21} and in the theory of electrical and communication networks \cite{20}, when it turned out that some computationally difficult problems can be solved on cacti in polynomial time. We refer the reader to paper \cite{5} for some aspects of domination in cactus graphs and to \cite{6} for some enumerative results on matchings and independent sets in chain cacti \cite{doslic, ARS}. 

A cactus graph is a connected graph in which no edge lies in more than one cycle. Consequently, each
block of a cactus graph is either an edge or a cycle. If all blocks of a cactus $G$ are cycles of the same size $m$,
the cactus is $m$-uniform.

 %As usual, we use $\lceil x \rceil$, $\lfloor x\rfloor$ for the smallest integer greater than or equal to $x$ and the largest integer less than or equal to $x$, respectively. In this article, we denote $\{1,2,\ldots,n\}$ simply by $[n]$.

\medskip

 The paper is organized as follows.
 In the next section,  we consider the number of independent dominating sets of triangular cactus. In Sections 3 and 4, we count the number of independent dominating sets of chain of squares  and chain of hexagonal cacti, respectively.

%%%%%%%%%%%%%%%%%%%%%%%%%%%%%%%%%%%%%%%%%%%%%%%%%%%%%%%%%%%%%%%%%%%%%%%%%%%%%%%%%
\section{Chain triangular cactus}

\begin{figure}[!ht]
\hspace{4cm}
\includegraphics[width=7.3cm,height=1.5cm]{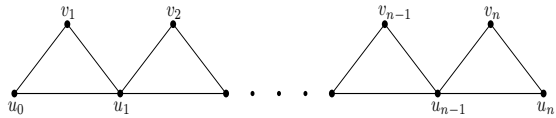}
\caption{ \label{cactus} The chain triangular cactus.}
\end{figure}

 In this section, we consider triangular cactus and call the number of triangles in $G$,  the length of the chain. An example of a chain triangular cactus is shown in Figure \ref{cactus}.
Obviously, all chain triangular cacti of the same length are isomorphic. Hence, we denote the chain triangular cactus of length $n$ by $T_n$. It is easy to see that for every natural numbers $n$, we have $\gamma_i(T_n)=\lfloor\frac{n+1}{2}\rfloor$.
 We investigate the generating function for the number of independent dominating sets of $T_n$.

\begin{theorem}
The generating function for the number of independent dominating sets of $T_n$ is given by
\[ T(x)=\frac{x+x^2}{1-x-x^2}.
\]
\end{theorem}
\begin{proof}
Let us consider $T_n$, labeled in the way  shown in Figure \ref{cactus}, and denote the number of independent dominating sets in $T_n$  by $t_n$. Each independent dominating set in $T_n$ either does or does not contain vertex $u_n$. Let us denote by $t^{\prime}_n$ the number of independent dominating sets that contain $u_n$, and by $t^{\prime\prime}_n$ the number of independent dominating sets that do not contain $u_n$. Hence, $t_n = t^{\prime}_n + t^{\prime\prime}_n$. 

Now we find recurrences for $t^\prime_n$ and $t^{\prime\prime}_n$. 

Obviously, Each independent dominating set in $T_n$ counted by $t^{\prime\prime}_n$ can be extended to an independent dominating set in $T_{n+1}$ counted by $t^{\prime}_{n+1}$ and we have $t^{\prime}_{n+1}=t^{\prime\prime}_{n}$.\\
Each independent dominating set in $T_n$ counted by $t^\prime_n$ is an independent dominating set in $T_{n+1}$ counted by $t^{\prime\prime}_{n+1}$. Further, an independent dominating set in $T_n$ counted by $t^{\prime\prime}_n$ can be extended to an independent dominating set in $T_{n+1}$ counted by $t^{\prime\prime}_{n+1}$ also in  only one way, by including $v_{n+1}$. Hence,
\[t^{\prime\prime}_{n+1}= t^{\prime}_{n}+t^{\prime\prime}_{n}.\]
We have obtained the system
\begin{eqnarray*}
t^{\prime\prime}_{n+1}&=& t^{\prime}_{n}+t^{\prime\prime}_{n}\\
t^{\prime}_{n+1}&=&t^{\prime\prime}_{n}
\end{eqnarray*} 
with the initial conditions $t^\prime_1=1$ and $t^{\prime\prime}_1=2$.

Now we introduce three generating functions, $T^\prime(x)=\sum_{n\geq 0} t^\prime_{n+1}x^n$ and $T^{\prime\prime}(x)= \sum_{n\geq 0} t^{\prime\prime}_{n+1}x^n$. By multiplying both equations in the above system through by $x^n$ and then summing over $n\geq 0$, the system can be translated into a linear system for two unknown generating functions: 
\[\begin{array}{lllll}
(1-x)T^{\prime\prime}(x)&-&xT^\prime(x)&=&1\\
T^\prime(x) &-& xT^{\prime\prime}(x)&=&0
\end{array}
\]
We obtain 
\[T^\prime(x)=\frac{x}{1-x-x^2},~~~T^{\prime\prime}(x)=\frac{1}{1-x-x^2}.\]
Finally, by adding $T^{\prime}(x)$ and $T^{\prime\prime}(x)$ and multiplying the sum by $x$ we obtain the generating function for the
sequence $t_n$, and $T(x)=\sum_{n\geq 0}t_nx^n$.\qed
\end{proof}

Since $T(x)$ is a rational function, we can conclude that the numbers $t_n$ satisfy a second order linear recurrence with constant coefficients. The initial conditions can be verified by direct computations. The following corollary gives the recurrence relation of $t_n$.
\begin{corollary}
For every $n\geq 3$, the number $t_n$ of independent dominating sets in $T_n$ is given by
\[ t_n= t_{n-1} +t_{n-2},
\] with the initial conditions $t_0=2$ and $t_1=3$.
\end{corollary}
The recurrence relation for the number $t_n$ of  independent dominating sets in $T_n$ is the Fibonacci recurrence, so an extremely good approximation to $t_n$ will be
\[t_n\sim \frac{1}{\sqrt{5}}(\frac{1+\sqrt{5}}{2})^n.\]
\section{Chains of squares}
In this section we count independent dominating sets in two classes of chains with all internal squares of the same
type. We start with the para-chain.

\subsection{Para-chain $Q_n$}

%*******************************************************
\begin{figure}[!ht]
\hspace{3.9cm}
\includegraphics[width=8.5cm,height=2cm]{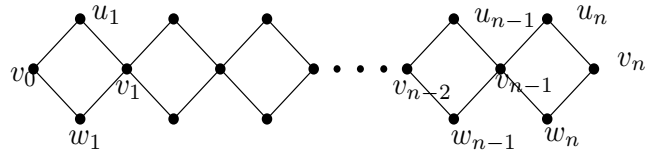}
\caption{ \label{pcactus} Labeled para-chain square cactus graphs. }
\end{figure}

\begin{theorem}\label{qn}
The generating function for the number of independent dominating sets of $Q_n$ is given by
\[ Q(x)=\frac{1+x^2}{(1-x)^2-x^3}.
\]
\end{theorem}
\begin{proof}
We consider a para-chain of length $n$, labeled as shown in Figure \ref{pcactus}. The number of independent dominating sets in $Q_n$ denoted by $q_n$. Each independent dominating set in $Q_n$ either does or does not contain vertex $v_n$. By $q^{\prime}_n$ we denote the number of independent dominating sets that contain $v_n$, and by $q^{\prime\prime}_n$ the number of independent dominating sets that do not contain $v_n$. Again, $q_n = q^{\prime}_n + q^{\prime\prime}_n$. In addition, we denote by $q^{\prime\prime\prime}_n$ the number of sets that are independent and not dominating in $Q_n$, but can be extended to an independent dominating set in $Q_{n+1}$. Clearly, such sets do not dominate $v_n$, and they must include $v_{n-1}$, since this vertex is necessary to dominate $u_n$ and $w_n$. Hence, they are counted by $q^\prime_{n-1}$ and we have $q^{\prime\prime\prime}_n=q^{\prime}_{n-1}$.

Now we find recurrences for $q^\prime_n$, $q^{\prime\prime}_n$ and $q^{\prime\prime\prime}_n$. 

Obviously, Each independent dominating set in $Q_n$ counted by $q^{\prime}_n$ and $q^{\prime\prime}_n$ can be extended to an independent dominating set in $Q_{n+1}$ counted by $q^{\prime}_{n+1}$  in  only one way. %Further, an independent dominating set in $Q_n$ counted by $q^{\prime\prime}_n$ can be extended to an independent dominating set in $Q_{n+1}$ counted by $q^{\prime}_{n+1}$ also in  only one way.
  We have the recurrence for $q^\prime_n$, 
\[q^\prime_{n+1}= q^\prime_n+q^{\prime\prime}_n=q_n.\]

Now we need a recurrence for $q^{\prime\prime}_n$. Each independent dominating set in $Q_n$ counted by $q^{\prime\prime}_n$ can be extended to an independent dominating set in $Q_{n+1}$ counted by $q^{\prime\prime}_{n+1}$ in  only one way, by   including both $u_{n+1}$ and $w_{n+1}$ and the same is valid for the sets counted by $q^{\prime\prime\prime}_n$. Hence,
\[q^{\prime\prime}_{n+1}= q^{\prime\prime}_{n}+q^{\prime\prime\prime}_{n}.\]
Finally, we have shown above that  $q^{\prime\prime\prime}_{n+1}=q^{\prime}_{n}$. We have obtained the system
\begin{eqnarray*}
q^\prime_{n+1}&=&q^\prime_n +q^{\prime\prime}_n\\
q^{\prime\prime}_{n+1}&=& q^{\prime\prime}_{n}+q^{\prime\prime\prime}_{n}\\
q^{\prime\prime\prime}_{n+1}&=&q^{\prime}_{n}
\end{eqnarray*} 
with the initial conditions $q^\prime_1=q^{\prime\prime}_1=q^{\prime\prime\prime}_1=1$.

We introduce three generating functions, $Q^\prime(x)=\sum_{n\geq 0} q^\prime_{n+1}x^n$, $Q^{\prime\prime}(x)= \sum_{n\geq 0} q^{\prime\prime}_{n+1}x^n$ and $Q^{\prime\prime\prime}(x)= \sum_{n\geq 0} q^{\prime\prime\prime}_{n+1}x^n$. By multiplying all equations in the above system through by $x^n$ and then summing over $n\geq 0$, the system can be translated into a linear system for three unknown generating functions:
\begin{eqnarray*}
(1-x)Q^{\prime}(x)&-&xQ^{\prime\prime}(x)=1\\
(1-x)Q^{\prime\prime}(x) &-& xQ^{\prime\prime\prime}(x)=1\\
Q^{\prime\prime\prime}(x) &-& xQ^{\prime}(x)=1
\end{eqnarray*}
We obtain
\[Q^\prime(x)=\frac{1+x^2}{(1-x)^2-x^3},~~Q^{\prime\prime}(x)=\frac{1}{(1-x)^2-x^3},~~Q^{\prime\prime\prime}(x)=\frac{1-x+x^2}{(1-x)^2-x^3}.\]
Finally, by adding $Q^{\prime}(x)$ and $Q^{\prime\prime}(x)$ and multiplying the sum by $x$ we obtain the generating function for the
sequence $q_n$ and $Q(x)=\sum_{n\geq 0}q_nx^n$.\qed
\end{proof}
%%%%%%%%%%%%%%%%%%%%%%%%%%%%%%%%%%%%%%%%%%%%%%%%%%%%%%

We can conclude that the numbers $q_n$ satisfy a third order linear recurrence with constant coefficients.  
The following corollary gives the recurrence relation of $q_n$.
\begin{corollary}
For every $n\geq 4$, the number $q_n$ of independent dominating sets in $Q_n$ is given by
\[ q_n= 2q_{n-1} -q_{n-2}+q_{n-3},
\] 
with the initial conditions $q_1=2,~ q_2=4$ and $q_3=7$.
\end{corollary}

\subsection{Ortho-chain $S_n$}

\begin{figure}[!ht]
\hspace{3.9cm}
\includegraphics[width=8.3cm,height=2.7cm]{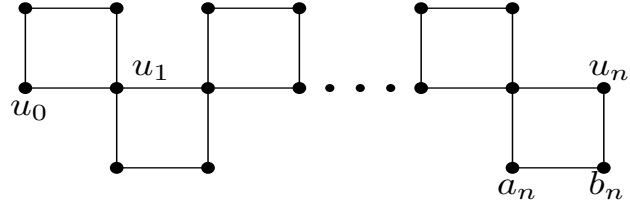}
\caption{ \label{scactus}  Labeled  ortho-chain square $S_n$. }
\end{figure}

%*******************************************************
Here we count the number of independent dominating sets in ortho-chain $S_n$ which  shown in Figure \ref{scactus}. 

\begin{theorem}\label{on}
The generating function for the number of independent dominating sets of $S_n$ is given by
\[ S(x)=\frac{1}{1-2x}.
\]
%which is valid for $|x|<\frac{1}{2}$.
\end{theorem}
\begin{proof}
We consider an ortho-chain of length $n$, labeled as shown in Figure \ref{scactus}. The number of independent dominating sets in $S_n$ denoted by $s_n$, and the number of independent dominating set in $S_n$ containing and not containing vertex $u_n$ are denoted by $s^{\prime}_n$ and  $s^{\prime\prime}_n$, respectively. Finally, we denote by $s^{\prime\prime\prime}_n$
the number of sets that are not independent dominating set in $S_n$, but can be extended to an independent dominating set in $S_{n+1}$. 

Analogously, before we obtain the system of recurrences for  $s^{\prime}_n$, $s^{\prime\prime}_n$ and $s^{\prime\prime\prime}_n$;
have the system
\begin{eqnarray*}
s^\prime_{n+1}&=&s^{\prime\prime}_n +s^{\prime\prime\prime}_n\\
s^{\prime\prime}_{n+1}&=& s^{\prime}_{n}+s^{\prime\prime}_{n}\\%=o_n\\
s^{\prime\prime\prime}_{n+1}&=&s^{\prime\prime}_{n}+s^{\prime\prime\prime}_{n} 
\end{eqnarray*} 
with the initial conditions $s^\prime_1=s^{\prime\prime}_1=1$.
We obtain $s_n=2s_{n-1}$ for $n\geq 1$ with the initial conditions $s_0 =1$, so  we have  the generating function $S(x) =\sum_{n\geq 0}2^nx^n$, and result follows.\qed
\end{proof}

\subsection{Examles}

\begin{figure}[!ht]
\hspace{3.9cm}
\includegraphics[width=8.5cm,height=4.5cm]{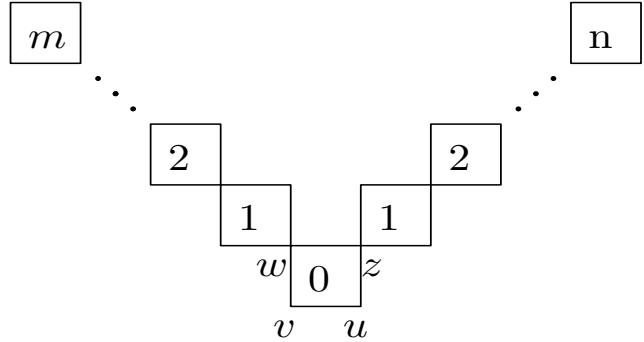}
\caption{ \label{dmn} A labeled para-chain with a single ortho-defect. }
\end{figure}

It would be interesting, for example, to show that the para-chain and the ortho-chain are extremal with respect to the number of independent dominating sets. As a step in that direction, we will investigate the effect of a single ortho-defect in a para-chain. The situation is shown in Figure \ref{dmn}. We denote this graph by $P_{mn}$ and the number of independent dominating sets in it by $p_{mn}$. 

%*******************************************************

In order to compute $p_{mn}$, we must consider two cases: whether an independent dominating set contains one of vertices $u$ and $v$. 
Let us look at the case when an independent dominating set in $P_{mn}$ contains $u$. There are $q_{m}(q^{\prime\prime}_n+q^{\prime\prime\prime}_n)$ such sets do not contain $z$. By symmetry, there must be $q_{n}(q^{\prime\prime}_m+q^{\prime\prime\prime}_m)$ independent dominating sets that contain $v$. 
 By the proof of theorem \ref{qn}, and summing the above contributions we obtain the following result: 
 \begin{eqnarray*}
 p_{mn}&=& q_{m}(q^{\prime\prime}_n+q^{\prime\prime\prime}_n)+q_{n} (q^{\prime\prime}_m+q^{\prime\prime\prime}_m)\\
 &=&q_{m}q^{\prime\prime}_{n+1}+q_{n}q^{\prime\prime}_{m+1}.
\end{eqnarray*} 
 
\begin{figure}[!ht]
\hspace{3.9cm}
\includegraphics[width=9.cm,height=2.cm]{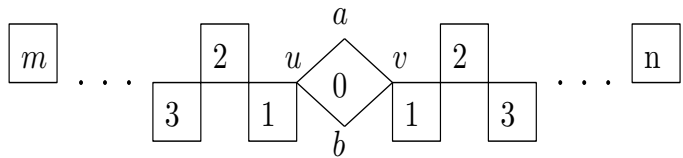}
\caption{ \label{omn} A labeled ortho-chain with a single para-defect. }
\end{figure}

%%%%%%%%%%%%%%%%%%%%%%%%%%%%%%%%%%%%%%%%%%%%%%%%%%%%%
Now, we  investigate the effect of a single para-defect in an ortho-chain. The situation is shown in Figure \ref{omn}. We denote this graph by $S_{mn}$ and the number of independent dominating sets in it by $s_{mn}$. 

%*******************************************************

In order to compute $s_{mn}$, we must consider three cases: whether an independent dominating set contains one or none of vertices $u$ and $v$. The case of none is the simplest: if an independent dominating set in $S_{mn}$ contains none of $u,~ v$,
then it must contain both of $a$ and $b$, and there are  $s^{\prime\prime}_n\cdot (s^{\prime\prime}_m+s^{\prime\prime\prime}_m)+ s^{\prime\prime\prime}_n\cdot (s^{\prime\prime\prime}_m+ s^{\prime\prime}_m)$ independent dominating sets. 

Let us now look at the case when an independent dominating set in $s_{mn}$ contains $u$ but not $v$. There are exactly $s^{\prime}_m\cdot s^{\prime\prime}_n$  such sets. By symmetry, there must be $s^{\prime\prime}_m\cdot s^{\prime}_n$ independent dominating sets that contain $v$ but not $u$. 
 By the proof of theorem \ref{on}, and summing the above contributions we obtain the following result: 
 \begin{eqnarray*}
 s_{mn}&=& (s^{\prime\prime}_n+ s^{\prime\prime\prime}_n)\cdot (s^{\prime\prime\prime}_m+ s^{\prime\prime}_m)+s^{\prime}_m \cdot s^{\prime\prime}_n+s^{\prime\prime}_m\cdot s^{\prime}_n\\
 &=&s_n\cdot s_m+2s_{m-1}\cdot s_{n-1}.
% &=& o^{\prime\prime}_m o_n+ o^{\prime\prime\prime}_no_m+ 2o^{\prime\prime}_no^{\prime}_m.
\end{eqnarray*}  
%%%%%%%%%%%%%%%%%%%%%%%%%%%%%%%%%%%%%%%%%%%%%%%%%%

%%%%%%%%%%%%%%%%%%%%%%%%%%%%%%%%%%%%%%%%%%%%%%%%%%%5
\section{Chain hexagonal cacti}

%*******************************************************
\begin{figure}[!ht]
\hspace{1.cm}
\includegraphics[width=12.8cm,height=4.8cm]{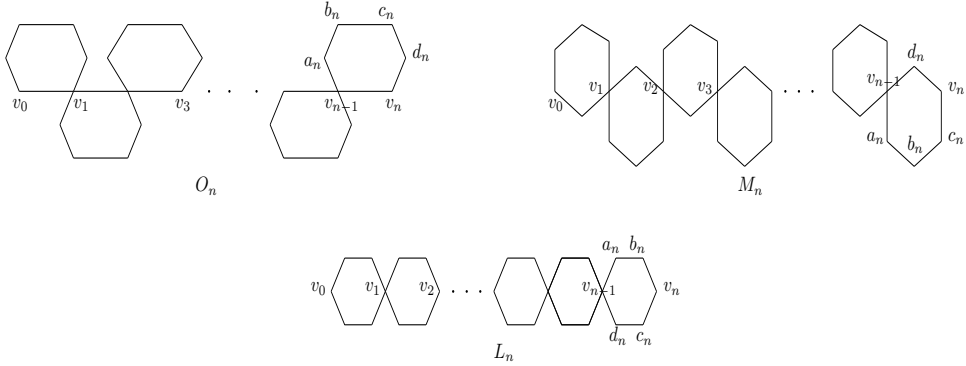}
\caption{ \label{hexagon} An Ortho-, meta-, and para-chain hexagonal cacti of length $n$. }
\end{figure}

In this section we investigate the number of independent dominating sets in three families of chain hexagonal cacti. 
The ortho-chain of length $n$ is denoted by $O_n$, and the meta-chain and the para-chain of length $n$ are denoted by $M_n$ and  $L_n$, respectively.  
Note that the independence domination number of Ortho-, and meta-chain hexagonal cacti of length $n$ is  $\lceil\frac{3n}{2}\rceil$.

We consider an ortho-chain hexagonal cacti of length $n$, labeled as shown in Figure \ref{hexagon}. The number of independent dominating sets in $O_n$ denoted by $o_n$, and the number of independent dominating set in $O_n$ containing and not containing vertex $v_n$ are denoted by $o^{\prime}_n$ and  $o^{\prime\prime}_n$. Finally, we denote by $o^{\prime\prime\prime}_n$
the number of sets that are not independent dominating set in $O_n$, but can be extended to an independent dominating set in $O_{n+1}$. 

Analogously, before we obtain the system of recurrences for  $o^{\prime}_n$, $o^{\prime\prime}_n$ and $o^{\prime\prime\prime}_n$;
have the system
\begin{eqnarray*}
o^\prime_{n+1}&=&2o^{\prime\prime}_n +2o^{\prime\prime\prime}_n\\
o^{\prime\prime}_{n+1}&=& 2o^{\prime}_{n}+2o^{\prime\prime}_{n}+o^{\prime\prime\prime}_{n}\\
o^{\prime\prime\prime}_{n+1}&=&o^{\prime\prime}_{n}+o^{\prime\prime\prime}_{n} 
\end{eqnarray*} 
with the initial conditions $o^\prime_1=2$ and $o^{\prime\prime}_1=3$.
We introduce the corresponding generating functions, $O^\prime(x)=\sum_{n\geq 0} o^\prime_{n+1}x^n$, $O^{\prime\prime}(x)= \sum_{n\geq 0} o^{\prime\prime}_{n+1}x^n$ and $O^{\prime\prime\prime}(x)= \sum_{n\geq 0} o^{\prime\prime\prime}_{n+1}x^n$ and obtain a linear system for them;
\[\begin{array}{lllllll}
O^{\prime}(x)&-&2xO^{\prime\prime}(x)&-&2xO^{\prime\prime\prime}(x)&=&2\\
(1-2x)O^{\prime\prime}(x) &-& 2 xO^{\prime}(x)&-&xO^{\prime\prime\prime}(x)&=&3\\
(1-x)O^{\prime\prime\prime}(x)  &-&xO^{\prime\prime}(x)&&&=&1.\\
\end{array}\]
We obtain \[O^\prime(x)=\frac{2+2x}{1-3x-3x^2},~~O^{\prime\prime}(x)=\frac{3+2x}{1-3x-3x^2},~~O^{\prime\prime\prime}(x)=\frac{{1+x}}{1-3x-3x^2}.\] 
Finally, by adding $O^{\prime}(x)$ and $O^{\prime\prime}(x)$ and multiplying the sum by $x$ we obtain the generating function for the sequence $o_n$ and $O(x)=\sum_{n\geq 0}o_nx^n$.

\begin{theorem}
The generating function for the number of independent dominating sets of $O_n$ is given by
\[ O(x)=\frac{1+2x+x^2}{1-3x-3x^2}.
\]
\end{theorem}
The following corollary gives the recurrence relation for  $o_n$.
\begin{corollary}
For every $n\geq 3$, the number of independent dominating sets in the chain hexagonal cacti $O_n$, i.e., is given by
\[ o_n= 3o_{n-1} +3o_{n-2},
\] 
with the initial conditions $o_1=5$ and $o_2=19$.
\end{corollary}
%%%%%%%%%%%%%%%%%%%%%%%%%%%%%%%%%%%%%%%%%%%%%%%%%%%%%%%%%%%%%%%

Now, we consider a meta-chain hexagonal cacti of length $n$, labeled as shown in Figure \ref{hexagon}. The number of independent dominating sets in $M_n$ denoted by $m_n$, and the number of independent dominating set in $M_n$ containing and not containing vertex $v_n$ are denoted by $m^{\prime}_n$ and  $m^{\prime\prime}_n$, respectively. Also, we denote by $m^{\prime\prime\prime}_n$
the number of sets that are not independent dominating set in $M_n$, but can be extended to an independent dominating set in $M_{n+1}$. Clearly, such sets do not dominate $v_n$, and they must include $v_{n-1}$, since this vertex is necessary to dominate $d_n$. Hence, they are counted by $m^\prime_{n-1}$ and we have $m^{\prime\prime\prime}_n=m^{\prime}_{n-1}$.

Analogously, before we obtain the system of recurrences for  $m^{\prime}_n$, $m^{\prime\prime}_n$ and $m^{\prime\prime\prime}_n$;
have the system
\begin{eqnarray*}
m^\prime_{n+1}&=&2m^{\prime\prime}_n + m^{\prime}_n +m^{\prime\prime\prime}_n\\
m^{\prime\prime}_{n+1}&=& m^{\prime}_{n}+2m^{\prime\prime}_{n}+2m^{\prime\prime\prime}_{n}\\
m^{\prime\prime\prime}_{n+1}&=&m^{\prime}_{n}
\end{eqnarray*} 
with the initial conditions $m^\prime_1=2$ and $m^{\prime\prime}_1=3$.
Again, we introduce the corresponding generating functions, $M^\prime(x)=\sum_{n\geq 0} m^\prime_{n+1}x^n$, $M^{\prime\prime}(x)= \sum_{n\geq 0} m^{\prime\prime}_{n+1}x^n$ and $M^{\prime\prime\prime}(x)= \sum_{n\geq 0} m^{\prime\prime\prime}_{n+1}x^n$ and obtain a linear system for them;
\[\begin{array}{lllllll}
(1-x-x^2)M^{\prime}(x)&-&2xM^{\prime\prime}(x)&=&1+x\\
(1-2x)M^{\prime\prime}(x) &-& (x+ 2x^2)M^{\prime}(x)&=&1+2x.
\end{array}\]
We obtain \[M^\prime(x)=\frac{1+x+2x^2}{1-3x-x^2-2x^3}, M^{\prime\prime}(x)=\frac{2x+1}{1-3x-x^2-2x^3},\! M^{\prime\prime\prime}(x)=\frac{1-2x}{1-3x-x^2-2x^3}.\] 
Finally, by adding $M^{\prime}(x)$ and $M^{\prime\prime}(x)$ and multiplying the sum by $x$ we obtain the generating function for the sequence $m_n$ and $M(x)=\sum_{n\geq 0}m_nx^n$.

\begin{theorem}
The generating function for the number of independent dominating sets of $M_n$ is given by
\[ M(x)=\frac{1-x+2x^2}{1-3x-x^2-2x^3}.
\]
\end{theorem}
The following corollary gives the recurrence relation of $m_n$.
\begin{corollary}
For every $n\geq 3$, the number $m_n$ of independent dominating sets in the chain hexagonal cacti $M_n$ is given by
\[ m_n= 3m_{n-1} +m_{n-2}+2m_{n-3},
\] 
with the initial conditions $m_0=1,~m_1=5$ and $m_2=19$.
\end{corollary}

%%%%%%%%%%%%%%%%%%%%%%%%%%%%%%%%%%%%

Now, we consider a para-chain hexagonal cacti of length $n$, labeled as shown in Figure \ref{hexagon}. The number of independent dominating sets in $L_n$ denoted by $l_n$, and the number of independent dominating set in $L_n$ containing and not containing vertex $v_n$ are denoted by $l^{\prime}_n$ and  $l^{\prime\prime}_n$. Finally, we denote by $l^{\prime\prime\prime}_n$
the number of sets that are not independent dominating set in $L_n$, but can be extended to an independent dominating set in $L_{n+1}$. 

Analogously, before we obtain the system of recurrences for  $l^{\prime}_n$, $l^{\prime\prime}_n$ and $l^{\prime\prime\prime}_n$;
have the system
\begin{eqnarray*}
l^\prime_{n+1}&=&l^{\prime\prime}_n + l^{\prime}_n +l^{\prime\prime\prime}_n\\
l^{\prime\prime}_{n+1}&=& l^{\prime}_{n}+3l^{\prime\prime}_{n}+2l^{\prime\prime\prime}_{n}\\
l^{\prime\prime\prime}_{n+1}&=&l^{\prime\prime}_{n}+l^{\prime\prime\prime}_n
\end{eqnarray*} 
with the initial conditions $l^\prime_1=2$ and $l^{\prime\prime}_1=3$.
Again, we introduce the corresponding generating functions, $L^\prime(x)=\sum_{n\geq 0} l^\prime_{n+1}x^n$, $L^{\prime\prime}(x)= \sum_{n\geq 0} l^{\prime\prime}_{n+1}x^n$ and $L^{\prime\prime\prime}(x)= \sum_{n\geq 0} l^{\prime\prime\prime}_{n+1}x^n$ and obtain a linear system for them;
\[\begin{array}{lllllll}
(1-x)L^{\prime}(x)&-&xL^{\prime\prime}(x)&-&xL^{\prime\prime\prime}&=&2\\
(1-3x)L^{\prime\prime}(x) &-& xL^{\prime}(x)&-&2xL^{\prime\prime\prime}&=&3\\
(1-x)L^{\prime\prime\prime}&-&xL^{\prime\prime}(x)&&&=&1.
\end{array}\]
We obtain% \[L^\prime(x)=\frac{x^2-4x+2}{1-5x+4x^2-x^3}, L^{\prime\prime}(x)=\frac{3-2x}{1-5x+4x^2-x^3},\! L^{\prime\prime\prime}(x)=\frac{x^2-x+1}{1-5x+4x^2-x^3}.\] 
\[L^\prime(x)=\frac{x^3-3x^2-4x+2}{1-6x+9x^2-6x^3+x^4},~~~ L^{\prime\prime}(x)=\frac{3-5x+4x^2-x^3}{1-6x+9x^2-6x^3+x^4},\] \[L^{\prime\prime\prime}(x)=\frac{2x^2-2x+1}{1-6x+9x^2-6x^3+x^4}.\] 
Finally, by adding $L^{\prime}(x)$ and $L^{\prime\prime}(x)$ and multiplying the sum by $x$ we obtain the generating function for the sequence $l_n$ and $L(x)=\sum_{n\geq 0}l_nx^n$.

\begin{theorem}
The generating function for the number of independent dominating sets of $L_n$ is given by
\[ L(x)=\frac{1-x-5x^3+x^4}{1-6x+9x^2-6x^3+x^4}.
\]
\end{theorem}
The following corollary gives the recurrence relation of $l_n$.
\begin{corollary}
For every $n\geq 4$, the number $l_n$ of independent dominating sets in the chain hexagonal cacti $L_n$ is given by
\[ l_n= 6l_{n-1} -9l_{n-2}+6l_{n-3}-l_{n-4},
\] 
with the initial conditions $l_0=4,l_1=5,l_2=19$ and $l_3=76$.
\end{corollary}

%%%%%%%%%%%%%%%%%%%%%%%%%%%%%%%%%%%%%%%%%%%%%%%%%%%%%%%%
\bigskip

%\noindent{\bfseries{Acknowledgements}}
%	The authors acknowledge the financial support from  Iran National Science Foundation (INSF),  Tehran and Yazd University research affairs (research project INSF-YAZD 96010014).

%%%%%%%%%%%%%%%%%%%%%%%%%%%%%%%%%%%%%%%%%%%%%%%%%%%%%%%%%%
%%%%%%%%%%%%%%%%%%%%%%%%%%%%%%%%%%%%%%%%%%%%%%%%%%%%%%%%%

\end{document}